\documentclass[11pt]{amsart}
\usepackage[utf8]{inputenc}
\usepackage{lipsum}
\usepackage{mwe}
\usepackage{amsfonts}
\usepackage{amssymb}
\usepackage{mathrsfs}
\usepackage{latexsym}
\usepackage{graphicx}
\usepackage{fancyhdr}
\usepackage[english]{babel}
\usepackage{amsmath}
\usepackage{amsthm}
\usepackage{enumerate}
\usepackage{color}
\usepackage{cases}
\usepackage[lofdepth,lotdepth]{subfig}
\usepackage{cite}
\usepackage{enumitem}
\usepackage{xcolor}
\usepackage{mathptmx}
\usepackage{float}
\usepackage{hyperref}

\newcommand{\R}{\mathbb{R}}

\newcommand{\di}{\text{div }}

\newcommand{\D}[1]{\text{D}#1}

\newtheorem{theorem}{Theorem}[section]

\theoremstyle{definition}
\newtheorem{definition}[theorem]{Definition}

\numberwithin{equation}{section}

\newcommand{\bphi}{\bar{\phi}}
\newcommand{\lbphi}{\bphi_{avg}}
\newcommand{\tphi}{\phi_\theta}
\newcommand{\bmu}{\bar{\mu}}
\newcommand{\tmu}{\mu_\theta}
\newcommand{\bint}{\int_{\partial\Omega}}
\newcommand{\bu}{\bar{u}}
\newcommand{\tu}{u_\theta}
\newcommand{\AZ}{A_0^{-1}(\bphi - \lbphi)}

\title[Convergence Properties of PINNs for the NSCH System]{Convergence Properties of PINNs for the Navier-Stokes-Cahn-Hilliard System}

\author[K. Buck]{Kevin Buck}
\address[K. Buck]{Indiana University, USA}
\email{{\tt kevbuck@iu.edu}}

\author[R. Temam]{Roger Temam}
\address[R. Temam]{Indiana University, USA}
\email{\tt temam@iu.edu}



\begin{document}

\begin{abstract}
Approximating solutions to differential equations using neural networks has become increasingly popular and shows significant promise. In this paper, we propose a simplified framework for analyzing the potential of neural networks to simulate differential equations based on the properties of the equations themselves. We apply this framework to the Cahn-Hilliard and Navier-Stokes-Cahn-Hilliard systems, presenting both theoretical analysis and practical implementations. We then conduct numerical experiments on toy problems to validate the framework's efficacy in accurately capturing the desired properties of these systems and numerically estimate relevant convergence properties.
\end{abstract}

\maketitle

\section{Introduction}

Approximating solutions to differential equations using neural networks is a rapidly advancing field, driven by increasing amounts of easily available computational power and advancements in the machine learning techniques themselves.  Traditional numerical methods, while effective, may struggle with nonlinear problems, with irregular solutions, and with high dimensional simulation.

In recent years, neural networks, particularly Physics-Informed Neural Networks (PINNs), have shown significant promise in providing efficient and accurate approximations to differential equations in great generality.  PINNs integrate the governing physical laws described by differential equations into the training of Neural Networks in place of (unsupervised learning) or in addition to (supervised learning) data. This allows these methods to take advantage of the structure of neural networks as well as the modeling strength given by differential equations.  Despite these advancements, a comprehensive body of work for analyzing the effectiveness of these neural networks a priori remains underdeveloped.

Several studies have demonstrated the potential of neural networks in solving differential equations \cite{heattransfer, highspeed, cardiac, radiative, poroelasticity, atomistic, fasttime, crustal, fluids8020046}.  While these examples demonstrate the successful application of neural networks the overall success of such attempts has been inconsistent, highlighting the need for more analytical results to assess the reliability of these techniques.  Some existing results surrounding the analysis of these networks can be found in \cite{residuals, HAL}.

The Cahn-Hilliard and Navier-Stokes-Cahn-Hilliard systems are of particular interest due to their widespread applications in modeling phase separation and fluid dynamics, respectively. The Cahn-Hilliard equation describes the process of phase separation in binary mixtures, while the Navier-Stokes equations model the behavior of fluid flow. Together they provide a complete picture of the phase separation of a mixed fluid in motion.  However, coupling these systems presents additional challenges that necessitate robust analytical tools to ensure accurate simulations.

This paper aims to perform this analysis by proposing a simplified framework for analyzing the potential of PINNs based on the intrinsic properties of the equations themselves and apply them to these systems. Our objectives are to:

\begin{enumerate}
   \item Present the framework we use to assess the effectiveness of PINNs and address a few points of interest.
   \item Apply this framework to the Cahn-Hilliard and Navier-Stokes-Cahn-Hilliard systems.
   \item Perform numerical experiments on a toy problem to validate the efficacy of the framework as well as numerically approximate some relevant constants.
\end{enumerate}

The structure of the paper is as follows.  In section 2.1 and 2.2 we describe the basics of Neural Networks and Physics-Informed Neural Networks largely as presented in \cite{RAISSI2019686, Convergence}.  After discussing some of the potential issues we hope to address in section 2.3, we describe the analytical properties of a problem necessary to avoid these issues in 2.4. We then present our approach in the context of the Cahn-Hilliard, then Navier-Stokes Cahn-Hilliard systems in section 3 and 4.  Finally in section 5, we look at toy problems numerics to see how the analysis behaves in practice.

\subsection{Description of the Neural Network}
A Neural Network aims to approximate an unknown function with a composition of nonlinear ``Activation'' functions with many parameters. The optimal values of these parameters are then found through the minimization of an objective, or ``Loss" function in the parameter space.

Precisely, a (Feedforward) Neural Network is a function $u$ with parameter matrices $A_i$ and vectors $b_i$ (denoted collectively as $\theta$) is defined as
\begin{equation}
    u(x) := A_n (\sigma( A_{n-1}(\sigma(\dots \sigma(A_1 x + b_1 ))) +b_{n-1}   ) +b_n
\end{equation}
where $\sigma$ is a chosen non-polynomial function often taken to be a sigmoid or rectified linear unit ($ReLU$) function, applied componentwise.  We usually denote $u$ as $u_\theta$ to emphasize the dependence on the parameters $\theta$.

The given Loss function, which depends (often implicitly) on the network parameters $\theta$, is then minimized with respect to $\theta$. It is common to take the Loss function to depend on $\theta$ stochastically, or minimize the function stochastically.  The most popular minimization algorithms are gradient based methods, particularly the quasi-newton method L-BFGS algorithm \cite{Liu1989} or the first-order, stochastic ADAM \cite{adam} algorithm.  The appropriate gradients are calculated via forward or backward mode Automatic Differentiation, which allows for efficient and accurate computation of the gradients (described in detail in \cite{autodiff}).

It is well understood that neural networks with an arbitrarily large amount of parameters are dense in many function spaces, including the space of continuous functions.  So, by minimizing this objective function we aim to ``train" the network to approximate some unknown target function in those spaces. 

\subsection{Physics-Informed Neural Networks}
To use the structure of the Neural Network to simulate the solutions to a differential equation we recall the ``Physics-Informed Neural Network" or PINN, introduced in \cite{RAISSI2019686}. In this structure we use the network to approximate the unknown solution function to a known differential equation. In this context the network has inputs x and t, and outputs the value of the solution function(s). 

This type of network is defined primarily by its Loss function, which is derived from the residual of the equation we hope to approximate. For example, consider a generic partial differential equation (PDE):

\begin{gather} \label{Introduction:GenericPDE}
\begin{cases}{}
    Au = f \text{ in }\Omega\times (0,T) \\
    u(0,x) = u_0(t) \text{ in }\Omega \\
    u(x,t) = g(x,t) \text{ on }\partial\Omega\times (0,T)
\end{cases}
\end{gather}
The corresponding analytic (or continuous) loss function is:
\begin{equation}
    L(\theta) = ||Au_\theta-f||_{\Omega\times(0,T)} + ||u_\theta(0,x) - u_0(x)||_\Omega + ||u-g||_{\partial\Omega\times (0,T)} 
\end{equation}
where $u_\theta$ is the network and $||\cdot||_V$ denotes an appropriate Sobolev or $L^p$ norm on the space $V$.  

In practice, this analytic loss function is replaced by an empirical (or discrete) loss function, which is a stochastic approximation of the analytic loss.  The empirical loss is computed by approximating the norms in the analytic loss using a large number of randomly chosen sample points, as in a Monte-Carlo approximation. The relationship between the empirical Loss function and the analytic Loss function is explored thoroughly in \cite{Convergence}.  This paper's results enable us to focus our analysis on the analytic loss function, as the analytic loss is bounded close to the empirical loss observed in practice.  

Additionally, regularization terms are often added to this Loss function to penalize the network if it deviates too far from the expected regularity of solutions.  This encourages the approximation to satisfy desired regularity requirements earlier in the minimization process, which is important for our later analysis. For example, if we expect the solution function to satisfy $L^2(0,T;H^1(\Omega))$ regularity, we use a regularized loss function given by

\begin{equation} \label{I:reg}
    L_{reg}(\theta) = L(\theta) + \lambda||u_\theta||_{L^2(0,T;H^1(\Omega))}
\end{equation}

\subsection{Analytical Concerns}
In order to prove and analyze the efficacy of these networks in approximating the solution to a differential equation, and to ease the implementation of these methods there are several concerns which need to be addressed.

A first concern is the viability of the Loss function described in the PINN framework above. In linear algebra, it is well-known that the minimization of a residual is not always an appropriate way to approximate the solution of a problem (See Appendix A). The setting of differential equations is much more complex, what reason do we have to believe that minimizing a residual is appropriate here?

A second concern involves the actual minimization procedure. While many numerical algorithms, such as L-BFGS and ADAM, are constantly improving for minimizing loss functions, there remains a need for analysis to explore under what conditions a PDE might produce a particularly stubborn local minimum or low gradient and importantly how to avoid them when they appear. Effective minimization algorithms can reduce the likelihood of worst-case scenarios during the optimization process. Some results on this phenomenon are discussed in \cite{krishnapriyan2021characterizing}.

In this paper, we focus solely on the first concern, leaving the second concern for future work. Our aim is to discuss the worst-case scenario of an optimization scenario and provide bounds for this worst case. First, we analyze the issue generally and explain the terminology used in our analysis. We then apply this analysis to the Cahn-Hilliard and Navier-Stokes Cahn-Hilliard models and verify these results with a few very simple numerical experiments.

\subsection{Analytical Framework}
For the described PINN framework, we study the relation between the value of the loss function and the actual error of the function approximation. The weakest positive result would be that as the value of the Loss function approaches zero, so does the actual error. We call this the property the consistency of a PINN.

\begin{definition}
    Let $u_\theta$ denote a Physics Informed Neural Network with parameters $\theta$ satisfying $Loss(\theta) < \epsilon$.  We say that the network is \textbf{consistent} in the $\Omega\times(0,T)$ norm if and only if
\begin{equation}
    ||u-u_\theta||_{\Omega\times(0,T)} \to 0 \text{ as } \epsilon \to 0.
\end{equation}
where $u$ is the exact solution to the underlying equation.
\end{definition}

This property is extremely minimal for sensible approximation of a differential equation by a Neural Network. In this way, this property can be viewed as similar to the consistency of a numerical scheme.  A PINN which does not satisfy this property has no hope of simulation the system, as even a process which minimizes the Loss to zero will not necessarily yield an accurate approximation if terminated early.

A more robust property is that of convergence:
\begin{definition}
    Let $u_\theta$ denote a Physics Informed Neural Network with parameters $\theta$ satisfying $Loss(\theta) < \epsilon$ for $\epsilon$ sufficiently small.  The network \textbf{converges with order $n$} in the $\Omega\times(0,T)$ norm if and only if
    \begin{equation}
    ||u-u_\theta||_{\Omega\times(0,T)} \leq C \epsilon^n + o(\epsilon^n)
\end{equation}
We call this constant $C$ the \textbf{Condition Number} of the PINN with respect to the order of convergence, and $\epsilon$ the radius of convergence.
\end{definition}

The details provided in this bound are useful in computational practice, as they can be used to estimate how small the Loss needs to be to achieve a desired level of accuracy. Importantly, it is not sufficient to study the value of $n$, as if $C$ is sufficiently large, achieving an accurate approximation through the minimization of the Loss is hopeless. For example, the results in \cite{Biswas2022} establish a promising rate of convergence for the Navier-Stokes system, however at the cost of a very large condition number.  Only by knowing the order of convergence along with the condition number can we determine an appropriately small value of the loss function for a given problem.  This type of Condition Number also arises in Linear Algebra, and is very important for accurate computation there. We describe this in detail in Appendix A.

The ambiguity surrounding the phrase ``for sufficiently small epsilon" should be acknowledged. For varying sizes of epsilon, we expect the sharpness of the bound to change (potentially drastically). We do not consider this problem in this paper as it is closely related to the problem of stubborn local minima in the neural network space, which is connected more closely to optimization methods (issue two of our analytical concerns).  Even weights yielding a neural network which is not a true local minima can have exceptionally small gradients, which in practice causes issues nearly as severe.

\subsubsection{Language of Operators}
It is important to note that the order of convergence and the condition number are properties of the differential equations (or more accurately the solution operator) rather than properties of the network itself. We can thus formulate these convergence properties in the language of operators. Denote the solution operator to (\ref{Introduction:GenericPDE}) by $S$, which takes inputs the forcing function $f$, the boundary conditions $g$, and the initial condition $u_0$. Note that what is regarded as an input and what is regarded as parameter can be changed depending on the context of the problem. Importantly, the PDE must satisfy existence properties for such an operator to exist, and uniqueness for the operator to be invertible. Let 
\begin{equation}
S(f, g, u_0) = u.
\end{equation}
Then let ($\delta f$, $\delta g$, $\delta u_0$) represent small perturbations in each of these input components and define $\delta u$ so that 
\begin{equation}
    S(f+\delta f, g + \delta g, u_0 + \delta u_0) = u + \delta u
\end{equation}  
Then the convergence property can be written as follows:

\begin{definition}
Let $u_\theta$ again denote a PINN with parameters $\theta$ satisfying $Loss(\theta)<\epsilon$ for $\epsilon$ sufficiently small.  The network converges with order in the $\Omega\times(0,T)$ norm if and only if
\begin{equation}
    ||\delta u|| \leq C \left(||(\delta f, \delta u_0, \delta g)||\right)^n + o\left(||(\delta f, \delta u_0, \delta g)||\right)^n
\end{equation}
\end{definition}

Here $C$ is the condition number and $n$ the rate of convergence, the same as above. This also clarifies that the value of the Condition number can depend on $f,g,$ and $u_0$.

\section{Cahn-Hilliard}
The Cahn-Hilliard equations describe the phase separation of a two-phase fluid through a diffuse interface model. This model captures the evolution of the concentration of two mixed substances through the evolution of the order parameter $\phi$, which represents the difference in the concentration of the two fluids. The interface between the fluids can then be determined by observing the level sets of the order parameter, thus recovering the sharp interface without tracking it directly.

Since $\phi$ is a difference of concentrations, we restrict to the case where $\phi(x,t)\in[-1,1]$. Solutions satisfying this condition are called ``physical solutions." It is well known that physical solutions to the Cahn-Hilliard model exist, so we restrict to this case for the remainder of the discussion.

The fourth order Cahn-Hilliard equations are represented as a system of two coupled second order equations with coupling term $\mu$:
\begin{gather}
    \begin{cases}{} \label{CH:MainEq}
        \partial_t\phi + (u\cdot\nabla)\phi = \Delta\mu &\text{ in }\Omega\\
        \mu = -\Delta\phi + \Psi'(\phi) &\text{ in }\Omega 
    \end{cases}
\end{gather}
with the Initial and Boundary Conditions:
\begin{gather}
\begin{cases}{}
    \phi(\cdot, 0) = \phi_0 &\text{ in } \Omega \label{CH:MainEq:IC} \\
    \frac{\partial\phi}{\partial n} = 0 \text{ , }  \frac{\partial\mu}{\partial n} =0 &\text{ on } \partial\Omega\times (0,T)
\end{cases}
\end{gather}
Note here that $u$, which here is taken as a parameter and denotes the velocity of the fluid, satisfies
\begin{align*}
    \nabla \cdot u = 0  \text{ in } \Omega\\
    u = 0 \text{ on } \partial\Omega
\end{align*}

We assume that $\Omega$ is a closed and bounded subset of $\R^2$. The decoupling $\mu$ denotes the chemical potential. The function $\Psi(s)$ denotes the free energy density. There are several forms of $\Psi$ which are of interest however for the Cahn-Hilliard section of this report we will use the Landau Potential, given by $\Psi(s) = \frac{1}{4}(s^2-1)^2$, for our model. For the analysis of Navier-Stokes Cahn-Hilliard, a more general function is employed which we describe there.

\subsection{The Cahn-Hilliard Physics-Informed Neural Network}

To approximate these equations, we define the network by its loss function, given by
\begin{equation} \label{Loss1}
\begin{split}
    L(\tphi) = &\alpha_1||\partial_t\tphi + (u\cdot\nabla)\tphi - \Delta\mu_\theta ||^2_{L^4(0,T;L^2(\Omega)} + \alpha_2||\tphi(\cdot, 0) - \phi_0 ||^2_{L^2(\Omega)} \\
    &+ \alpha_3||\frac{\partial\tphi}{\partial n}||^2_{L^4(0,T;L^2(\partial\Omega))} + \alpha_4||\frac{\partial\tmu}{\partial n}||^2_{L^4(0,T;L^2(\partial\Omega))} \\
    &+ \lambda_1 || \tphi||_{L^2(0,T;H^1(\Omega))} + \lambda_2 ||\tmu||_{L^2(0,T;H^1(\Omega))}
\end{split}
\end{equation}
where $\tmu = -\Delta\tphi - \Psi'(\tphi) $, and $\alpha_i$ are constant parameters which determine the relative learning rate of each term. The loss is a function of only the parameters $\theta$. Here the spatial and temporal derivatives of $\tphi$ are calculated through automatic differentiation, as described in the introduction.  We note that the final two terms of the Loss are not strictly necessary but enforce a regularity mirroring that of a true solution, as described by \ref{I:reg}.  We require $H^1$ regularization of the network so that we may apply trace theorems in our analysis.  Even if not enforced by the loss, this regularity will occur naturally through the use of a smooth activation function in the network.

For a network simulating the Cahn-Hilliard equation, the output is a single value in the interval $[-1,1]$, which approximates the value of the solution function $\phi$ at the point $(x,t)$.  Implementing this restriction is straightforward in practice, and can be done for instance by using a \textit{tanh} activation function on the final layer of the network.

It is important to discuss the consistent difficulties of simulating the Cahn-Hilliard (and similar equations such as Allen-Cahn) with PINNs.  These simulations were carried out using a special Runga-Kutta Neural Network in \cite{RAISSI2019686} at least partially because of the failure of a typical PINN structure to reliably learn the correct solution.  This difficulty arises because of the particularly stubborn local behavior of these operators in particular.  Since the difficulty is in the training itself, it is necessary to attempt a more sophisticated training method to alleviate the problem.  This came in the form of time-adaptive training \cite{wight2020solving} and similarly motivated space-adaptive training \cite{haitsiukevich2023improved} which solve these particular training issues for a large class of problems. Other effective methods, such as those in \cite{MATTEY2022114474}, have been developed to assist with the training of Cahn-Hilliard and Allen-Cahn equations.   Despite these advancements, poorly conditioned problems exist and remain to be solved.  However, selecting an appropriate training method can significantly reduce the negative effects of a poorly conditioned problem.  Another recent result with similar goals is found in \cite{zhang2024priori}, where other a priori estimates for the convergence of a Neural Network approximating Cahn-Hilliard are derived.

\subsection{Results for Cahn-Hilliard}
The main result for the Cahn-Hillard system is found in the following theorem:
\begin{theorem}
    Let $\tphi(x,t)$ be the neural network with outputs in [-1, 1] and $\phi(x,t)$ a weak, physical solution to the Cahn-Hilliard equations (2). If the Loss associated with the network, defined by (\ref{Loss1}), satisfies 
    \begin{equation} \label{Theorem1:LossBound}
    L(\tphi)<\epsilon^2,
    \end{equation}
    then
    \begin{equation} \label{Theorem1:FinalBound}
        ||\tphi(t)-\phi(t)||_{L^2(\Omega)}^2 \leq\epsilon^2e^{C_\lambda t}\left(\frac{1}{\alpha_2} + \left(\frac{1-e^{-2C_\lambda t}}{2} \right) \left( \frac{\lambda^2}{2\alpha_1} + \frac{C_9}{\alpha_3} + \frac{C_9}{\alpha_4} \right)\right).
    \end{equation}
    Where $\lambda > 0$ can be chosen arbitrarily. $C_\lambda$ depends only on $\lambda$ and $\Omega$. Additionally $C_\lambda\searrow C$ as $\lambda \nearrow\infty$ and $C_\lambda\nearrow \infty$ as $\lambda \searrow 0$.
\end{theorem}
This theorem guarantees that a network approximating the solution converges at least with order 1, and describes the condition number up to constants relating to the domain.  The proof of this theorem is found in Appendix B.

\section{Navier-Stokes-Cahn-Hilliard}
It is natural to couple the Cahn-Hilliard equations with the Navier-Stokes equations to replace the parameter function $u$ with a simulated fluid flow. We use the coupling known as ``Model H," named for Hohenberg and Halperin. This model can be derived both through Continuum Mechanics and through an Energy Variational approach and is given by

\begin{gather} 
    \begin{cases}{} \label{NSCH:MainEq}
        u_t + (u\cdot\nabla)u -\nabla\cdot(\nu(\phi)Du) + \nabla P = \mu\nabla\phi \\  
        \nabla\cdot u = 0 \\
        \phi_t + (u\cdot\nabla)\phi = \nabla\cdot\left(\nu(\phi)\nabla\mu\right) \\
        \mu = -\Delta\phi + \Psi'(\phi)
    \end{cases}
\end{gather}

with initial and boundary conditions:

\begin{gather}
    \begin{cases}{}
    u(\cdot,0) = u_0 \text{ , } \phi(\cdot, 0) = \phi_0 &\text{ in } \Omega \label{NSCH:MainEq:ICBC} \\
    u=0 \text{ , } \frac{\partial\phi}{\partial n} = 0 \text{ , }  \frac{\partial\mu}{\partial n} =0 &\text{ on } \partial\Omega\times (0,T)
    \end{cases}
\end{gather}

As before $\Omega\subseteq \R^2$.  Here, $u(x,t):\Omega\times(0,T)\to\R^2$ denotes the velocity, and $\phi(x,t):\Omega\times(0,T)\to\R$ represents the order parameter. The quantity $\nu$ denotes the kinematic viscosity, $P$ is the pressure. Additionally, $D$ denotes the symmetric gradient operator.

\subsection{Navier-Stokes Cahn-Hillard PINN}
We aim to prove convergence for a PINN approximating these equations.  To do this, we define the network with inputs $x$ and $t$ and outputs $\tu, \tphi, P_\theta$ with the restrictions that $\text{div }\tu = 0$ and
\begin{equation}
        \tu = 0\text{,  }\partial_n\tmu = \partial_n\tphi = 0 \text{ on }\partial\Omega\times (0,T).
\end{equation}
It is reasonable to construct a network which meets boundary conditions (and even the initial condition) exactly as shown in  \cite{Sun2019SurrogateMF, chudomelka2020deep, ko2022convergence}.  The restriction that the divergence condition be met exactly is also attainable by the methods shown in \cite{RichterPowell2022NeuralCL}.

With these restrictions, we define the Loss of this Neural Network by 
\begin{equation}
\begin{split}
    L(\tu, \tphi, P_\theta) :=& ||\partial_t\tu + (\tu\cdot\nabla)\tu - \text{div }(\eta(\tphi)D\tu) + \nabla P_\theta - \tmu\nabla\tphi||_{L^4(0,T;L^2(\Omega))} \\
    &+ ||\partial_t\tphi + \tu\nabla\tphi - \Delta\tmu ||_{L^4(0,T;L^2(\Omega))} + ||\tu(\cdot, 0) - u_0||_{L^2(\Omega)} \\
    &+ ||\tphi(\cdot, 0) - \phi_0 ||_{L^2(\Omega)} + \lambda_1 ||\tu||_{L^2(0,T;V_\sigma)} + \lambda_2 ||\tphi||_{L^2(0,T;H^2(\Omega))}
\end{split}
\end{equation}
where $\tmu = -\Delta\tphi + \Psi'(\tphi)$ and the function spaces are defined as in Appendix B.

We note similarly to the Cahn-Hilliard case that the final two terms of the Loss are not strictly necessary, but instead mirror the known regularity (as described in Appendix D) of weak solution by the process described in (\ref{I:reg}). 

There are some existing results for the Navier-Stokes equations.  In \cite{Biswas2022}, the authors derive a bound similar to the ones we describe here for the Navier-Stokes equation.  In this analysis the divergence free condition is removed, making the analysis more complete.  This result also suggests how removing the divergence free condition would impact the rate of convergence.  It should be noted that the derived condition number (not described as such in the reference) is extremely large, in addition to the decreased rate of convergence.  However, this type of small convergence rate should be expected for such a general class of complex problems.

Additionally, an in depth discussion of some other practical concerns of training a network for approximating the Navier-Stokes equations can be found \cite{deryck2023error}.  Here there is a detailed discussion of the number of nodes required for the network, as well as the expected decrease of the Loss.

Many papers attempt to simulate the Navier Stokes equation with various twists and on different problems.  A general but detailed overview can be found in \cite{Jin_2021}.  Many additional references on Navier Stokes can be found in the survey paper \cite{Cuomo2022}.

For the Navier-Stokes Cahn-Hillard system itself, less has been done.  Some numerical results and discussions can be found \cite{doi:10.1137/18M1223459, GIORGINI2020194, NSCHattractors, Chen2022}.

\subsection{Results for Navier-Stokes-Cahn Hilliard}

\begin{theorem}
    Let $(u, \phi, P, \mu)$ be a weak solution to (\ref{NSCH:MainEq})-(\ref{NSCH:MainEq:ICBC}). Let $(\tu, \tphi, P_\theta)$ be an output of the Neural Network described above satisfying
    \begin{equation}
        L(\tu, \tphi, P_\theta) < \epsilon,
    \end{equation}
    and $\tphi(x,t) \in [-1, 1].$  Then $||u - \tu||_\# + ||\phi - \tphi ||_* \to 0$ as $\epsilon\to 0$.
\end{theorem}

This theorem shows that the Neural network is consistent in the to be described ($H^{-1}$-like) norms.  The proof of this theorem can be found in Appendix C.

\section{Numerical Results}
We perform some numerical experiments on a toy problem Navier-Stokes-Cahn-Hilliard system using a PINN in order to numerically estimate the rate of convergence and condition number.  The numerics for Cahn-Hillard have been performed by others in papers such as \cite{RAISSI2019686, wight2020solving, zhang2024priori}, though full understanding of how to efficiently simulate this system is still very much ongoing.

We use a Loss function slightly different from the one used in analysis. Firstly we do not include the regularizer. Including the regularizer will only serve to increase the performance. Secondly, we do not use a network which is divergence free.  We also split the learning of $\phi$ and $\mu$, though the analysis is not affected by this change.  Finally, we use a large coefficient in front of the Initial Condition terms.  This leads to the following analytic loss function.

\begin{equation}
\begin{split}
    L(\tu, \tphi, P_\theta) :=& ||\partial_t\tu + (\tu\cdot\nabla)\tu - \text{div }(\eta(\tphi)D\tu) + \nabla P_\theta - \tmu\nabla\tphi||_{L^2(0,T;L^2(\Omega))} \\
    &+ ||\tmu + \Delta\tphi - \Psi'(\tphi)||_{L^2(0,T;L^2(\Omega))} \\
    &+ ||\partial_t\tphi + \tu\nabla\tphi - \Delta\tmu ||_{L^2(0,T;L^2(\Omega))}  + ||\nabla\cdot\tu||_{L^2(0,T;L^2(\partial\Omega))}\\
    &+ 1000(||\tu(\cdot, 0) - u_0||_{L^2(\Omega)}+ ||\tphi(\cdot, 0) - \phi_0 ||_{L^2(\Omega)})    \\
    &+ ||\frac{\partial\tphi}{\partial n} ||_{L^2(0,T;L^2(\partial\Omega))} + ||\frac{\partial\tphi}{\partial n} ||_{L^2(0,T;L^2(\Omega))}
\end{split}
\end{equation}

To train the network we use the time-adaptive method described in \cite{wight2020solving} in conjunction with a learning rate scheduler. We use ADAM optimization at each step. This leads to efficient minimization of the chosen Loss.  We observe the loss on an arbitrarily chosen problem, different from the problem we perform our next experiments on.  For this problem we cannot compute the analytic solution exactly, so we cannot measure the exact error.  However, we can see that the Loss decreases to reasonable levels in Figure 3 (around $3 * 10^{-3}$ neglecting the 1000 coefficient on the initial condition).

\begin{figure}[H]
    \centering
    \includegraphics[height=3in]{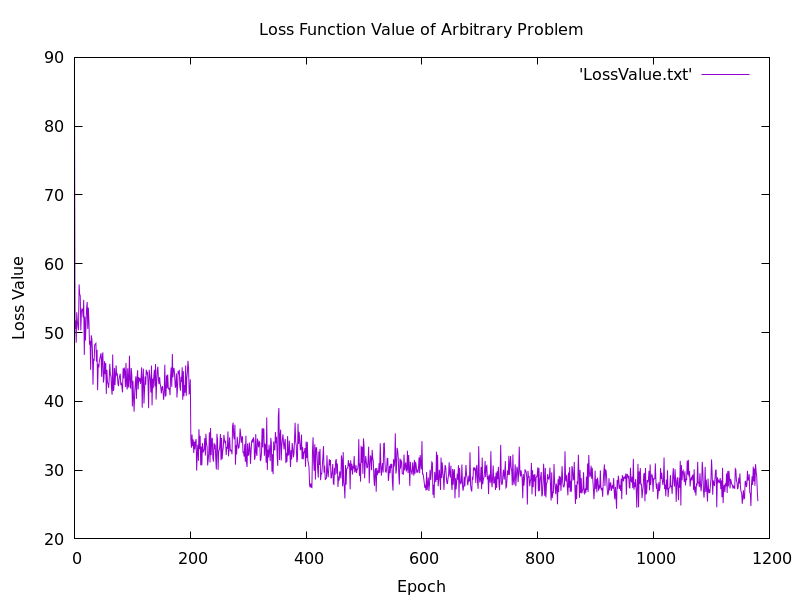}

    \caption{We see the loss decreasing on an arbitrarily chosen problem different from our other experiments.  For this problem, the Loss (without the 1000 weight on the initial condition) is on the order of $3*10^{-3}$}
\end{figure}\label{fig:decreasingLoss}

For further experiments, we need to effectively measure the error.  In order to accomplish this we choose a problem for which the exact analytic solution can be computed.  Specifically, we choose forcing terms so that the exact solution is given by:
\begin{center}
\begin{equation}
\begin{split}
    u(x,t) &= \begin{bmatrix}
        -x_2 t \\ x_1 t 
    \end{bmatrix} \\
    \phi(x,t) &= \sin(x_1t) + \sin(x_2t)
\end{split}
\end{equation}
\end{center}

This solution is simple, but allows us to study the relationship between the actual error and the value of the networks loss function at least in one particular case.  We then train many different networks on this problem, terminating the minimization at different points to study the relationship between loss and error at relatively high and low values of the loss.  We simulate until final time $T=2$.  This leads to the results in Figure 2, which plot the weighted loss value against the actual $L^2(0,T;L^2(\Omega))$ (computed stochastically through Monte-Carlo) error for many different training attempts of different length.

\begin{figure}[H]
    \centering
    \includegraphics[height=3in]{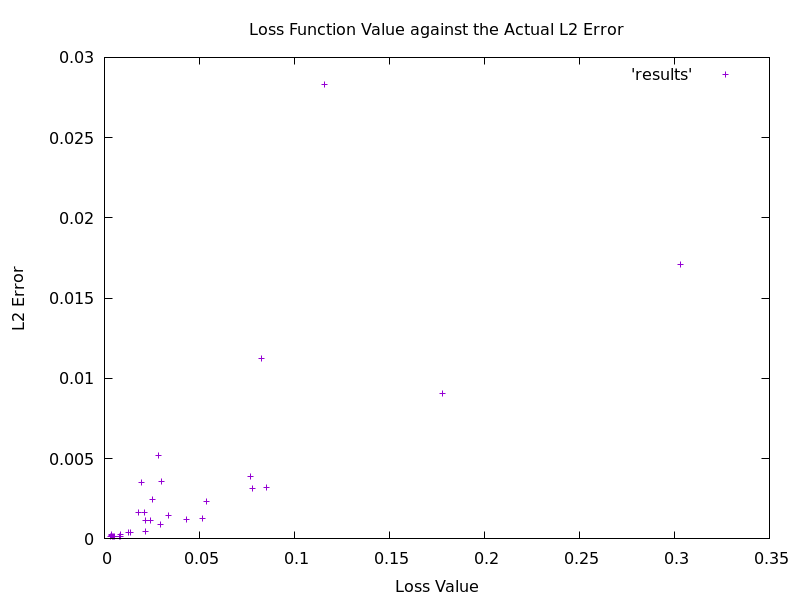}

    \caption{The $L^2$ error at the final timestep plotted against the Emperical Loss Value for the described PINN modeling the NSCH system given above, before any processing of data.  We can see that the actual $L^2$ error of the simulation is quite small for small Loss values, implying the desired convergence.  In some cases the error is as small as about $10^{-4}$}
\end{figure}\label{fig:preprocessed}

The results allow us to numerically estimate the rate of convergence and condition number.  By plotting the same data on a logscale in Figure 3, we can observe that the line of best fit (shown also in Figure 3) is given by
\begin{equation}
\log(\text{L2 error}) = 1.153 \log\text{(Loss function value)} - 1.044.
\end{equation}
This suggests a rate of convergence of $1.153$ and a condition number of about $.09$ for this particular problem.  This is quite good, much better than would be expected of a generic problem.  However, this is due purely to the fact that we use such a simple toy problem.

\begin{figure}[H]
    \centering
    \includegraphics[height=3in]{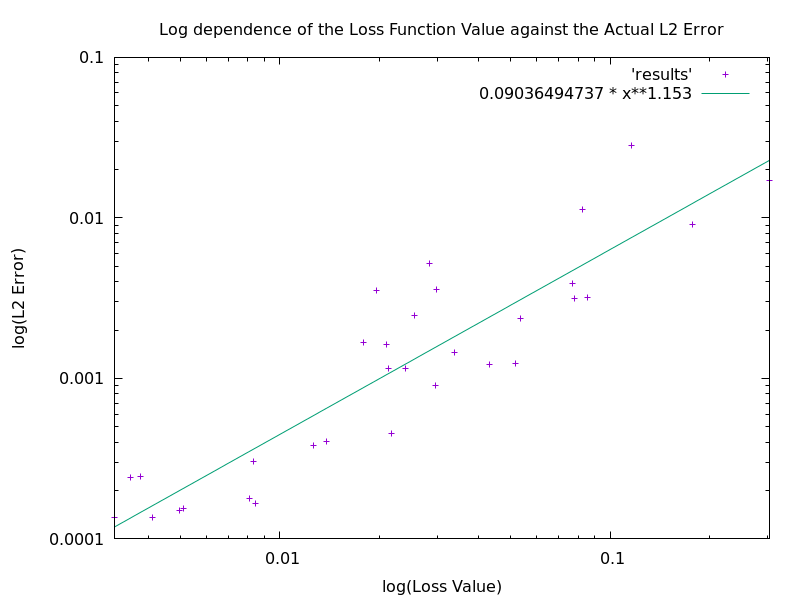}

    \caption{The $L^2$ error vs the Emperical Loss Value for the described PINN modeling the NSCH system, after the processing of data to a base 10 logscale.  We can see that the data descends roughly in parallel with the derived line of best fit with condition number $.09$ and rate of convergence $1.153$}
\end{figure}\label{fig:postprocess}

All of the code used can be found at the following GitHub: \url{https://github.com/kevmbuck/Analysis-of-NSCH-PINNs}

\section{Conclusion and Further Thoughts}
In this study we have proposed a simplified framework for analyzing the convergence properties of Physics Informed Neural Networks, and applied that framework to the Cahn-Hilliard and Cahn-Hilliard Navier-Stokes systems.  For the Cahn-Hilliard system we obtained a convergence result as well as estimates on the condition number of the system.  In the case of the full Navier-Stokes Cahn-Hilliard, we obtain a consistency result and test numerically the rate of convergence for a simple problem.

It is our hope that this framework can be applied to obtain useful information about the effectiveness of a PINN in simulating particular problems.  This is especially helpful in interpreting the results of an attempted PINN simulation, by giving concrete bounds on the accuracy given the value of the Loss function at the final optimization step.

A weakness of the above analysis is the several simplifying assumptions made in the particular analysis of the Navier-Stokes Cahn-Hilliard system, as well as the lack of truly comprehensive numerics.  Furthermore we neglect in our analysis the numerical optimization process itself.  By assuming a worst-case optimization scenario, our results are very broad.  However, it could be that in practice a particularly effective optimization procedure can drastically outperform our analysis and thus limit its usefulness.

In the future we hope to explore results which bound the condition numbers and rates of convergence of a system for particular classes of initial data or other parameter regimes of the system.  By having a firm understanding of how the condition number and rate of convergence are affected by system parameters, we can further our understanding of the strengths and weaknesses of PINN simulation and when or how it should be applied.  This idea could also be expanded from our analysis of a worst case to include an expected case with a given optimization procedure.  This would allow us a much more thorough understanding of which problems are appropriate for various optimization methods.

We additionally want in the future to get more specific convergence results for simpler systems of PDE.  For these simpler systems, it is possible that we could achieve more precise estimates for the condition number and rates of convergence dependent on qualities of the initial condition, domain, and other system parameters.  We would also like to be able to implement and test these methods on a more comprehensive set of benchmark problems, and see if the rate of convergence and condition numbers can be reliably determined by numerics.

\section*{Appendix A: Linear Algebra}
Here we explore one example which illustrates the potential failure of using the minimization of the residual to solve problems in Linear Algebra.

We aim to approximate the solution to the simple linear equation
\begin{equation}
A\vec{x} = \vec{b}
\end{equation}
through the minimization of the residual, so the problem becomes
\begin{equation}
  \min_{\vec{y}} ||A\vec{y}-b||
\end{equation}

Since we do not minimize the residual exactly, the residual of our approximate solution $\vec{y}$ has some positive value.  We thus define $\delta b :=  A\vec{y} - b$ which is nonzero.  Then we decompose $\vec{y}$ into two components $\vec{x}$, the true solution, and $\delta x$, the deviation from the true solution.  The value of the residual of our approximate solution $\vec{y}$ is thus equal to $||\delta b||$.  This allows us to define the Condition Number of the matrix $A$ to be the smallest number $M$ such that 
\begin{equation}
    \frac{||\delta x||}{||\vec{x}||} \leq M \frac{||\delta b||}{||\vec{b}||}, \forall b\neq 0, \delta b > 0
\end{equation}.

If the condition number is sufficiently large, it is possible that a small residual can lead to an extremely inaccurate approximation.  This can happen even for problems which might otherwise seem unproblematic.  For example, consider the case where
\begin{equation}
A = 
\begin{bmatrix}
10& 7 &8 &7 \\
7 &5 &6 &5  \\
8 &6& 10& 9  \\
7 &5 &9 &10 
\end{bmatrix}, \hspace{1cm}\vec{b}=
\begin{bmatrix}
    32 \\23\\ 33\\ 31
\end{bmatrix},\hspace{1cm} \vec{x} = \begin{bmatrix}
    1 \\ 1 \\ 1 \\ 1
\end{bmatrix}.
\end{equation}
Here $\vec{x}$ is the exact solution to the equation.  Additionally, $A$ has many desirable qualities.  It is invertible with inverse containing only integer entries (the largest of which has absolute value 68), and has determinant 1.  There is also no obvious reason a priori to believe that $\vec{b}$ should cause any difficulties.

There are however $\vec{y}$ with very small residual values which are not at all close to the true solution $\vec{x}$.  Take for example the $\vec{y}$ given below, then we have 

\begin{equation}
    \vec{y} = \begin{bmatrix}
        10.2 \\ -11.6 \\ 5.5 \\ 5.1
    \end{bmatrix}, \text{ which yields }
    \delta b = \begin{bmatrix}
        10^{-1} \\ -10^{-1} \\ 10^{-1} \\ -10^{-1}
    \end{bmatrix}, \text{ however }
    \delta x = \begin{bmatrix}
        9.2 \\ -12.6 \\ 4.5 \\ 4.1
    \end{bmatrix}.
\end{equation}

In this example although $\vec{y}$ has very small residual (in the sup norm equal to only $10^{-1})$, the actual error is very large.  This type of problem necessarily has a large Condition number, and is thus called ill-conditioned.  For this reason, it is not appropriate to solve the problem through residual minimization.

\section*{Appendix B: Function Spaces and Notation}
We define $\Omega\subset\R^2$ a bounded domain with smooth boundary.  We denote $W^{k,p}(\Omega)$ the Sobolev spaces of functions with $k$ distributional derivatives in $L^p(\Omega)$ and its norm $||\cdot||_{W^{k,p}(\Omega)}$.  We denote the Hilbert space $W^{k,2}(\Omega)$ by $H^k(\Omega)$ with norm $||\cdot||_{H^k}$.  Denote by $H_0^1(\Omega)$ the closure of $C^\infty_0(\Omega)$ with respect to the $H^1(\Omega)$ norm.  We denote it's dual $H^{-1}$.

We define $H=L^2(\Omega)$ with norm $||\cdot||$ and inner product $(\cdot, \cdot)$.  We denote $V=H^1(\Omega)$ with norm $||\cdot||_V$ and denote its dual $V'$.  Let $\langle\cdot,\cdot\rangle$ denote the duality product between $V$ and $V'$.  We denote by $u_{avg}$ the average value of $u$ on $\Omega$, given by $u_{avg} = |\Omega|^{-1}(u,1)$.

We denote by $C^\infty_{0,\sigma}(\Omega)$ the space of divergence free vector fields in $C^\infty_0(\Omega)$.  We define $H_\sigma$ the closer of $C^\infty_{0,\sigma}(\Omega)$ with respect to the $H$ norm, and $V_\sigma$ the closer of $C^\infty_{0,\sigma}(\Omega)$ with respect to the $V$ norm.  We use the same notation for inner product and norm on these space as for $H$ and $V$ respectively.

We note that $V_\sigma$ has the inner product and norm $(u,v)_{V_\sigma} = (\nabla u, \nabla v)$ and $||u||_{V_\sigma} = ||\nabla u||$.  We denote the dual of $V_\sigma$ by $V_\sigma'$.  By Korn's Inequality, the norm $u\to ||D u ||$ is an equivalent norm on $V_\sigma$, where $D$ again here denotes the symmetric gradient.

\section*{Appendix C: Proof of Cahn-Hilliard Results}
\subsection*{C.1 Preliminaries}

%\subsubsection*{C.1.1 Regularity of Solutions}
We assume that the initial datum is appropriately smooth to allow the solution functions $\tmu$, $\tphi$ to be in $H^1(\Omega)$ in space.  This allows us to apply the Trace theorem.

\subsection*{C.2 Proof of Theorem}
\begin{proof}
First we recognize that by defining
\begin{align}
    &f_1(x,t) := \partial_t\tphi + (u\cdot\nabla)\tphi - \Delta\mu_\theta \\
    &f_2(x) := \tphi(\cdot,0) - \phi_0 \\
    &f_3(x) := \frac{\partial\tphi}{\partial n}   \\
    &f_4(x) := \frac{\partial\tmu}{\partial n}
\end{align}
$\tphi$ satisfied the following equation.
\begin{equation}
    \begin{cases}
        \partial_t\tphi + (u \cdot \nabla)\tphi = \Delta\mu_\theta + f_1 &\text{ in } \Omega\times (0,T) \\
        \mu_\theta = -\Delta\tphi + \Psi'(\tphi)
    \end{cases}
\end{equation}
with Initial and Boundary Conditions
\begin{equation}
    \begin{cases}
        \tphi(\cdot, 0) = \phi_0 + f_2 \text{ in } \Omega \\
        \frac{\partial\tphi}{\partial n} = f_3 \text{ , }  \frac{\partial\mu_\theta}{\partial n} = f_4 \text{ on } \partial\Omega\times (0,T)
    \end{cases}
\end{equation}

We proceed by analyzing the differences $\bar{\phi} = \phi-\tphi$ and $\bar{\mu} = \mu - \mu_\theta$.  Which satisfy a third system
\begin{equation} \label{MR:DiffEq}
    \begin{cases}
        \partial_t\bar{\phi} + (u \cdot \nabla)\bar{\phi} = \Delta\bar{\mu} - f_1 &\text{ in } \Omega\times (0,T) \\
        \bar{\mu} = -\Delta\bar{\phi} + \Psi'(\phi) - \Psi'(\tphi) &\text{ in } \Omega\times (0,T) \\
    \end{cases}
\end{equation}
with Initial and Boundary Conditions
\begin{equation} \label{bhi:BC}
    \begin{cases}
        \bphi(\cdot, 0) = f_2 \text{ in } \Omega \\
        \frac{\partial\bphi}{\partial n} = 0 \text{ , }  \frac{\partial\bmu}{\partial n} = 0 \text{ on } \partial\Omega\times (0,T)
    \end{cases}
\end{equation}
Now we attempt to estimate $\bar{\phi}$.  First, condense (\ref{MR:DiffEq}) into one equation
\begin{equation}
    \partial_t\bar{\phi} + (u\cdot\nabla) + \Delta^2\bar{\phi} = \Delta(\Psi'(\phi) - \Psi'(\tphi))-f_1.
\end{equation}
Now we multiple both sides by $\bar{\phi}$ and integrate over $\Omega$ to see
\begin{equation}
\begin{split}
    \int_\Omega \partial_t\bphi \cdot \bphi dx &+ \int_\Omega (u \cdot\nabla)\bphi \cdot \bphi dx+ \int_\Omega \Delta^2\bphi \cdot \bphi dx  \\
    &= \int_\Omega (\Delta(\Psi'(\phi)-\Psi'(\tphi))\cdot \bphi dx - \int_\Omega f\bphi dx.
\end{split}
\end{equation}
We label the terms in this equation as $I_1 + I_2 + I_3 = I_4 - (f_1, \bphi)_{L^2(\Omega)}$.  We handle each term individually.  First,
\begin{equation}
    I_1 = \int_\Omega \partial_t\bphi \cdot \bphi dx = \frac{1}{2} \frac{d}{dt}|\bphi|^2_{L^2(\Omega)}.
\end{equation}
Then the second term,
\begin{align*}
    I_2 = \int_\Omega (u \cdot\nabla)\bphi \cdot \bphi dx &= \sum_i \int_\Omega u_i \partial_i \bphi \cdot \bphi dx \\
    &= \sum_i \int_\Omega \frac{u_i}{2} \partial_i(\bphi^2) dx\\
    &= \sum_i \int_{\partial\Omega} u_i \cdot n_i \frac{\bphi^2}{2} dx- \sum_i \int_\Omega \partial_i u_i \frac{\phi^2}{2} dx \\
    &= 0
\end{align*}
since $u\cdot n = 0$ on $\partial\Omega$ and $\nabla\cdot u$ = 0 in $\Omega$.  

We handle the third term by integrating by parts twice and using the boundary condition
\begin{align}
\begin{split}
    I_3 &= \int_\Omega \Delta^2\bphi \cdot \bphi dx \\
    &= - \int_\Omega \nabla\bphi \cdot \nabla\Delta\bphi dx + \int_{\partial\Omega}\frac{\partial\Delta\bphi}{\partial n}\bphi d\Gamma \\
    &= \int_\Omega |\Delta\bphi |^2 dx + \int_{\partial\Omega}\frac{\partial\Delta\bphi}{\partial n}\bphi d\Gamma - \int_{\partial\Omega} \frac{\partial \bphi}{\partial n} \cdot \Delta\bphi \\
\end{split}
\end{align}
Breaking this up further term by term we have by noting
\begin{equation} \label{Theorem1:PsiBC}
    \frac{\partial}{\partial n} \Psi'(\phi) = 3\phi \frac{\partial\phi}{\partial n} - \frac{\partial\phi}{\partial n} = 0 \text{ on } \partial\Omega dx.
\end{equation}
and using equation (\ref{MR:DiffEq}) that
\begin{align}
\begin{split}
    \int_{\partial\Omega} \frac{\partial\Delta\bphi}{\partial n} \bphi &= \bint - \frac{\partial\bmu}{\partial n}\bphi + \frac{\partial}{\partial n}(\Psi'(\phi) - \Psi'(\tphi))\bphi \\
    &\geq -\left|\frac{\partial\bmu}{\partial n}\right|_{L^2(\partial\Omega)}|\bphi|_{L^2(\partial\Omega)} + \bint (3\tphi^2-1)\bphi \frac{\partial\tphi}{\partial n} \\
    &\geq -\left|\frac{\partial\bmu}{\partial n}\right|_{L^2(\partial\Omega)}|\bphi|_{L^2(\partial\Omega)} - C_1 |\bphi|_{L^2\partial\Omega}\left|\frac{\partial\tphi}{\partial n}\right|_{L^2(\partial\Omega)} \\
    &= -|\bphi|_{L^2(\partial\Omega)} \left( \left| \frac{\partial\bmu}{\partial n}\right|_{L^2(\partial\Omega)} + C_1 \left|\frac{\partial\tphi}{\partial n}\right|_{L^2(\partial\Omega)}\right) \\
    &\geq -C_{Tr}^2 |\bphi|_{L^2(\Omega)}^2 - \left( \left| \frac{\partial\bmu}{\partial n}\right|^2_{L^2(\partial\Omega)} + C_1 \left|\frac{\partial\tphi}{\partial n}\right|^2_{L^2(\partial\Omega)}\right).
\end{split}
\end{align}
Here we assume $\bphi\in H^1(\Omega)$, then apply the trace theorem with trace constant $C_{Tr}$.  Next we observe that
\begin{equation}
    -\bint \frac{\partial\bphi}{\partial n} \Delta\bphi \geq -\left|\frac{\partial\bphi}{\partial n} \right|_{L^2(\partial\Omega)}|\Delta\bphi|_{L^2(\partial\Omega)}. \\
\end{equation}
Assuming that $\bmu\in H^1(\Omega)$ we again apply the trace theorem, this time on $\Delta\bphi$ with trace constant $C_{Tr2}$, we see
\begin{align}
\begin{split}
    -\bint \frac{\partial\bphi}{\partial n} \Delta\bphi &\geq -C_{Tr2} \left|\frac{\partial\bphi}{\partial n} \right|_{L^2(\partial\Omega)}|\Delta\bphi|_{L^2(\Omega)} \\
    &\geq -C_2 \left|\frac{\partial\bphi}{\partial n} \right|^2_{L^2(\partial\Omega)} - \frac{1}{10} |\Delta\bphi|^2_{L^2(\Omega)}.
\end{split}
\end{align}
So together we get the following bound on $I_3$
\begin{equation}
    I_3 \geq \frac{9}{10}|\Delta\bphi|_{L^2(\Omega)}^2   -C_{Tr}^2 |\bphi|_{L^2(\Omega)}^2 - \left( \left| \frac{\partial\bmu}{\partial n}\right|^2_{L^2(\partial\Omega)}  + C_1 \left|\frac{\partial\tphi}{\partial n}\right|^2_{L^2(\partial\Omega)}+ C_2 \left|\frac{\partial\bphi}{\partial n} \right|^2_{L^2(\partial\Omega)}\right).
\end{equation}

Now for $I_4$,
\begin{align}
\begin{split}
    I_4 &= \int_\Omega (\Delta(\Psi'(\phi)-\Psi'(\tphi))\bphi dx \\
    &= \int_{\partial\Omega} \frac{\partial}{\partial n} (\Psi'(\phi)-\Psi'(\tphi))\bphi d\Gamma - \int_{\partial\Omega} \left(\Psi'(\phi) - \Psi'(\tphi)\right) \frac{\partial\bphi}{\partial n} d\Gamma \\
    &\hspace{2cm}+ \int_\Omega \left( \Psi(\phi) - \Psi(\tphi)\right)\Delta\bphi dx. \\
\end{split}
\end{align}
Here again we work term by term,
\begin{align}
\begin{split}
    \bint \frac{\partial}{\partial n} (\Psi'(\phi)-\Psi'(\tphi))\bphi d\Gamma &= \bint \bphi\frac{\partial}{\partial n} \Psi'(\tphi) \\
    &= \bint\bphi\left(3\tphi^2-1\right)\frac{\partial\tphi}{\partial n} \\
    &\leq C_3 \left|\frac{\partial\tphi}{\partial n}\right|_{L^2(\partial\Omega)}|\bphi|_{L^2(\partial\Omega)} \\
    &\leq C_3^2\left|\frac{\partial\tphi}{\partial n}\right|^2_{L^2(\partial\Omega)} + C^2_{Tr}|\bphi|_{L^2(\Omega)},
\end{split}
\end{align}
which uses the assumed boundedness of $\tphi$.  Then the second boundary term becomes
\begin{align}
\begin{split}
    - \bint \left(\Psi'(\phi) - \Psi'(\tphi)\right) \frac{\partial\bphi}{\partial n} d\Gamma &\leq \left|\frac{\partial\bphi}{\partial n}\right|_{L^2(\partial\Omega)}|\Psi'(\phi) - \Psi'(\tphi)|_{L^2(\partial\Omega)} \\
    &\leq \left|\frac{\partial\bphi}{\partial n}\right|^2_{L^2(\partial\Omega)} + |\Psi'(\phi) - \Psi'(\tphi)|^2_{L^2(\partial\Omega)}.
\end{split}
\end{align}
Now we analyze $|\Psi(\phi) - \Psi'(\tphi)|^2_{L^2(\partial\Omega)}$.
\begin{align}
\begin{split}
    |\Psi(\phi) - \Psi'(\tphi)|^2_{L^2(\partial\Omega)} &= \bint (\phi^3 - \phi - \tphi^3 + \tphi)^2 dx \\
    &= \bint \left(( \bphi^3 - \bphi + 3\tphi^2\bphi + 3\tphi\bphi^2\right)^2 dx \\
    &= \bint \left(( 3\tphi^2 + 3\tphi\bphi + \bphi^2 - 1)\bphi\right)^2 dx \\
    &\leq C_4\bint \bphi^2  dx \\
    &\leq C_4C_{Tr}|\bphi|^2_{L^2(\Omega)}
\end{split}
\end{align}
where $C_4 = \sup_{x\in\partial\Omega} 3\tphi^2 + 3\tphi\bphi + \bphi^2 - 1$.  Similarly, we get
\begin{equation}
    |\Psi'(\phi) - \Psi'(\tphi)|^2_{L^2(\Omega)} \leq C_5|\phi|^2_{L^2(\Omega)}
\end{equation}
for $C_5 = \sup_{x\in\Omega} 3\tphi^2 + 3\tphi\bphi + \bphi^2 - 1$ which we will use later.  Both these $C$ exist and are finite.

The final term for $I_4$ is bounded by
\begin{align}
    \int_\Omega \left( \Psi'(\phi) - \Psi(\tphi)\right)\Delta\bphi dx \leq \frac{1}{10}|\Delta\bphi|^2_{L^2(\Omega)}+ C_6  |\Psi(\phi) - \Psi'(\tphi)|^2_{L^2(\Omega)}.
\end{align}
So we get the final bound on $I_4$:
\begin{equation}
    I_4\leq  C_7|\phi|^2_{L^2(\Omega)} + \frac{1}{10}|\Delta\bphi|^2_{L^2(\Omega)} + C_3^2\left|\frac{\partial\tphi}{\partial n}\right|^2_{L^2(\partial\Omega)} + \left|\frac{\partial\bphi}{\partial n}\right|^2_{L^2(\partial\Omega)},
\end{equation}
where $C_7 = C_5(C_6+1)+C_{Tr}^2$.
Finally we need to bound $I_5$.  We use the simple bound
\begin{equation}
    I_5\leq |\bphi|^2_{L^2(\Omega)} + |f_1|^2_{L^2(\Omega)}.
\end{equation}

Combining the previous inequalities we get the following differential inequality:
\begin{equation}
\begin{split}
    \frac{1}{2} \frac{d}{dt}|\bphi|^2_{L^2(\Omega)} + \frac{8}{10}|\Delta\bphi|_{L^2(\Omega)} 
    - \left( \left| \frac{\partial\bmu}{\partial n}\right|^2_{L^2(\partial\Omega)}  + C_1 \left|\frac{\partial\tphi}{\partial n}\right|^2_{L^2(\partial\Omega)}+ C_2 \left|\frac{\partial\bphi}{\partial n} \right|^2_{L^2(\partial\Omega)}\right) \\
    \leq (C_7+C_{Tr}^2)|\phi|^2_{L^2(\Omega)}  + C_3^2\left|\frac{\partial\tphi}{\partial n}\right|^2_{L^2(\partial\Omega)} 
    + \left|\frac{\partial\bphi}{\partial n}\right|^2_{L^2(\partial\Omega)} -  (f_1,\bphi)_{L^2(\Omega)}.
\end{split}
\end{equation}
Then by simplifying and noting that $\frac{\partial\tphi}{\partial n} = -\frac{\partial\bphi}{\partial n}$, we see
\begin{equation}
\begin{split}
    \frac{1}{2} \frac{d}{dt}|\bphi|^2_{L^2(\Omega)}  \leq C_8|\phi|^2_{L^2(\Omega)}  + C_9\left( \left| \frac{\partial\bmu}{\partial n}\right|^2_{L^2(\partial\Omega)}  +  \left|\frac{\partial\bphi}{\partial n} \right|^2_{L^2(\partial\Omega)} \right) - (f_1,\bphi)_{L^2(\Omega)}.
\end{split}
\end{equation}

By Young's inequality, for any $\lambda>0$, we get
\begin{equation}
    \frac{d}{dt} |\bphi|^2_{L^2(\Omega)} \leq \left(C_8 + \frac{1}{2\lambda^2}\right)|\phi|^2_{L^2(\Omega)} +C_9\left( \left| \frac{\partial\bmu}{\partial n}\right|^2_{L^2(\partial\Omega)}  + \left|\frac{\partial\bphi}{\partial n} \right|^2_{L^2(\partial\Omega)} \right) +\frac{\lambda^2}{2}|f_1|^2_{L^2(\Omega)} .
\end{equation}
Note that $C_8 + \frac{1}{2\lambda^2} \leq C_\lambda$ with $C_\lambda$ dependent only on $\lambda$ and $\Omega$.  Then we apply Gronwall's Inequality
\begin{align} \begin{split}
    |\bphi(t)|^2_{L^2(\Omega)}& \leq |\bphi(0)|^2_{L^2(\Omega)}e^{C_\lambda t} + e^{C_\lambda t} \int_0^t e^{-C_\lambda s}\Bigg( \frac{\lambda^2}{2}|f_1(s)|^2_{L^2(\Omega)} \\
    &+ C_9\left| \frac{\partial\bmu}{\partial n}\right|^2_{L^2(\partial\Omega)}  + C_9 \left|\frac{\partial\bphi}{\partial n} \right|^2_{L^2(\partial\Omega)} \Bigg) ds.
\end{split} \end{align}
Finally we use the assumption (\ref{Theorem1:LossBound}) and Holder's inequality to see
\begin{align}
\begin{split}
    |\bphi(t)|^2_{L^2(\Omega)} &\leq |\bphi(0)|^2_{L^2(\Omega)}e^{C_\lambda t} + e^{C_\lambda t}\left(\int_0^t \left(e^{-C_\lambda s}\right)^2 ds\right)^\frac{1}{2} \Bigg(\frac{\lambda^2}{2}|f_1|^2_{L^4(0,T;L^2(\Omega))} \\
    & + C_9\left| \frac{\partial\bmu}{\partial n}\right|^2_{L^4(0,T;L^2(\partial\Omega))}+ C_9\left|\frac{\partial\bphi}{\partial n} \right|^2_{L^4(0,T;L^2(\partial\Omega))}  \Bigg)\\
    &\leq\frac{\epsilon^2}{\alpha_2}e^{C_\lambda t} + e^{C_\lambda t}\left(\frac{1-e^{-2C_\lambda t}}{2} \right) \left( \frac{\lambda^2\epsilon^2}{2\alpha_1} + \frac{C_9\epsilon^2}{\alpha_3} + \frac{C_9\epsilon^2}{\alpha_4} \right)\\
    &= \epsilon^2e^{C_\lambda t}\left(\frac{1}{\alpha_2} + \left(\frac{1-e^{-2C_\lambda t}}{2} \right) \left( \frac{\lambda^2}{2\alpha_1} + \frac{C_9}{\alpha_3} + \frac{C_9}{\alpha_4} \right)\right)
\end{split}
\end{align}
which finishes the result.
\end{proof}

\section*{Appendix D: Proof of NSCH}
\subsection*{D.1 Preliminaries}
\subsubsection*{D.1.1 Existence and Regularity of Solutions}
From \cite{doi:10.1137/18M1223459} we have the following theorems for existence and uniqueness of strong solutions to (\ref{NSCH:MainEq}).  First we define the notion of weak solution.

\begin{definition} \label{AppD:Reg}
    Let $T > 0$. Given $u_0\in H_\sigma$, $\phi_0\in V\cap L^\infty(\Omega)$ with $||\phi_0||\leq 1$ and $|\phi_{0, avg}| < 1$, a pair $(u, \phi)$ is a weak solution to $(\ref{NSCH:MainEq})-(\ref{NSCH:MainEq:ICBC})$ on $[0, T]$ if
    \begin{align} \begin{split}
        &u \in L^\infty (0, T; H_\sigma)\cap L^2(0, T; V_\sigma )\text{, } \partial_t u \in L^{\frac{4}{d}}(0, T, V_\sigma'), \\
        &\phi\in L^\infty (0,T;V)\cap L^2(0, T; H^2(\Omega))\cap H^1(0,T;V'), \\
        &\phi\in L^\infty (\Omega\times (0, T)), \text{ with } |\phi(x, t)| < 1 \text{ a.e. } (x,t)\in\Omega\times (0, T),
    \end{split} \end{align}
    and satisfies
    \begin{equation}  \begin{split}
            \langle \partial_t u, v\rangle + b(u, u, v) + (v(\phi) Du, Dv ) = (\mu\nabla\phi, v) \hfill \forall v\in V_\sigma \\
            \langle \partial_t\phi, v\rangle + (u\cdot\nabla\phi, v) + (\nabla\mu, \nabla v) = 0  \hfill \forall v\in V
    \end{split}  \end{equation}
\end{definition}

Then the following theorem, also from \cite{doi:10.1137/18M1223459}, guarantees the regularity necessary for analysis.

\begin{theorem}
    Assume that $u_0\in H_\sigma$, $\phi_0\in V\cap L^\infty(\Omega)$ with $||\phi_0||_{L^\infty(\Omega)} \leq 1$ and $|\phi_{0, avg}| < 1$.  Then, for any $T > 0$, there exists a weak solution $(u,\phi)$ to $(\ref{NSCH:MainEq}) - (\ref{NSCH:MainEq:ICBC})$ on $[ 0,T ]$ in the sense of the above definition such that
    \begin{align} \begin{split}
        u&\in C([0,T], H_\sigma), \\
        \phi&\in C([0,T], V)\cap L^4([0,T]; H^2(\Omega))\cap L^2([0,T];W^{2,p}(\Omega),
    \end{split} \end{align}
    where $2\leq p < \infty$.
\end{theorem}
This theorem guarantees the regularity of the weak solution necessary for our analysis.
\subsubsection*{D.1.2 Parameter Functions}
We take the kinematic viscosity to be given by
\begin{equation}
    \nu(z) = \nu_1 \left(\frac{1+z}{2}\right) + \nu_2 \left(\frac{1-z}{2}\right).
\end{equation}
We also require
\begin{equation}
    0 < \nu_* \leq \nu(z) \leq \nu^*.
\end{equation}
We take a more general version of $\Psi$ than in the Cahn-Hilliard case, instead taking any $\Psi\in C([-1,1])\cap C^3(-1,1)$ which is of the form
\begin{equation} \label{PF:alpha0}
    \Psi(z) = F(z) - \frac{\alpha_0}{2}z^2 \hspace{1cm}\forall z\in [-1,1]
\end{equation}
with
\begin{equation} \label{PF:alpha1}
    \lim_{z\to-1}F'(z) = -\infty, \hspace{1cm} \lim_{z\to 1}F'(z) = \infty, \hspace{1cm} F''(z)\geq \alpha_1 > 0
\end{equation}
It is additionally essential that 
\begin{equation} \label{PF:alpha}
\alpha_0-\alpha_1 = \alpha > 0.
\end{equation}

We define $F(z) := \infty$ for any $z\not\in[-1,1]$ and additionally require that $F$ be convex.  These assumptions can be thought of as minimal for a viable potential.

\subsubsection*{D.1.3 Inequalities}
Recall the Gagliardo-Nirenberg Inequality for $d=2$,
\begin{equation} \label{Ineq:GN}
    ||u||_{L^4(\Omega)} \leq C ||u||^{\frac{1}{2}}||u||_V^\frac{1}{2} \hspace{1cm} \forall u\in V
\end{equation}
and the Agmon Inequalities, again for $d=2$,
\begin{equation} \label{Ineq:AG}
    ||u||_{L^\infty(\Omega)}\leq ||u||^\frac{1}{2}||u||^\frac{1}{2}_{H^2(\Omega)} \hspace{1cm} \forall u\in H^2(\Omega)
\end{equation}
\begin{equation} \label{Ineq:BG}
    ||\nabla u ||_{L^4(\Omega)} \leq C||u||^\frac{1}{2}_{L^\infty(\Omega)} ||u ||^\frac{1}{2}_{H^2(\Omega)}.
\end{equation}
Finally, recall the Brezis-Gallout inequality,
\begin{equation}
    ||u||_{L^\infty(\Omega)} \leq C ||u||_V \left[ \log\left( c + \frac{||u||_{H^2(\Omega)}}{||u||_V}   \right) \right]^\frac{1}{2}
\end{equation}

\subsubsection*{D.1.4 Neumann Problems}
For $\lambda\geq 0$, consider the Neumann problem
\begin{gather}
    \begin{cases}{}
        -\Delta u + \lambda u = f \hspace{1cm} \text{in }\Omega \\
        \partial_n u = 0
    \end{cases}
\end{gather}
We introduce $B_\lambda \in L(V,V')$ defined by
\begin{equation}
    \langle B_\lambda u, v\rangle = \int_\Omega \nabla u \cdot \nabla v + \lambda uv dx \hspace{1cm}\forall u,v\in V
\end{equation}
Consider
\begin{equation}
    V_0 = \{ v\in V : v_{avg} = 0\}\text{,   }V_0' = \{f\in V' : f_{avg}=0 \}
\end{equation}
and recall that $V=V_0\oplus\R$, $V' = V_0'\oplus\R$.  Denote by $A_0$ the restriction of $B_0$ to an isomorphism from $V_0$ onto $V_0'$, denote $A_0^{-1}: V_0'\to V_0$ its inverse.  Since $\forall f\in V_0'$, $A_0^{-1}f$ is the unique $u\in V_0$ such that $\langle A_0 u, v\rangle = \langle f, v\rangle$ $\forall v\in V$.  Observe that
\begin{equation}
    \langle f, A_0^{-1}g\rangle = \int_\Omega \nabla (A_0^{-1})\cdot \nabla (A_0^{-1}g) dx \hspace{1cm} \forall f,g\in V_0'.
\end{equation}
Thus we can define
\begin{equation}
    ||f||_* = ||\nabla A_0^{-1}f|| = \langle f, A_0^{-1}f\rangle^\frac{1}{2}
\end{equation}
which is a norm on $V_0'$ which is equivalent to the natural one.  Furthermore for any $u\in H^1(0,T;V_0')$ we have
\begin{equation}
    \langle u_t(t), A_0^{-1}(u(t)) \rangle = \frac{1}{2} \frac{d}{dt} ||u(t)||_*^2 \text{ for almost every } t\in(0,T).
\end{equation}
Due to regularity we also have
\begin{equation}
    ||f||_* = ||\nabla A_0^{-1} f||_V \leq C||f|| \hspace{1cm}\forall f\in H\cap V_0'.
\end{equation}

\subsubsection*{D.1.5 Stokes Operators}
Consider the homogenous Stokes Problem
\begin{align}
    \begin{cases}{}
        -\Delta u + \nabla p = f \hspace{1cm} \text{ in } \Omega \\
        \di u = 0 \hspace{1.9cm}\text{ in }\Omega \\
        u = 0 \hspace{2.45cm}\text{ on }\partial\Omega
    \end{cases}
\end{align}

Let $A: V_\sigma\to V_\sigma'$ be such that
\begin{equation}
    \langle Au, v\rangle = (\nabla u, \nabla v) \hspace{1cm}\forall u,v \in V_\sigma
\end{equation}
This is an isomorphism, denote by $A^{-1}: V_\sigma'\to V_\sigma$ its inverse.  So given $f\in V_\sigma'$, $\exists!u=A^{-1}f\in V_\sigma$ such that
\begin{equation}
    (\nabla A^{-1}f, \nabla v) = \langle f, v\rangle \hspace{1cm} v\in V_\sigma.
\end{equation}
Thus $||f||_\# = ||\nabla A^{-1} f|| = \langle f, A^{-1}f\rangle^\frac{1}{2}$ is a norm on $V_\sigma$ equivalent to the natural one, which satisfies $\langle f_t(t), A^{-1}f(t)\rangle = \frac{1}{2}\frac{d}{dt} ||f(t)||^2_\#$ for almost every $t\in (0,T)$ for all $f\in H^1(0,T; V_\sigma')$.

To recover $p$, the well known result of De Rham implies $\exists p\in H$ such that $p_{avg} = 0$ and $\nabla p = \Delta u + f$ if $f\in H^{-1}(\Omega)$.  Additionally,
\begin{equation} \label{Ineq:pbound}
    ||p||\leq C||f||_{H^{-1}(\Omega)}
\end{equation}
for some $C>0$.

Assuming $f\in H$, there exists a unique $u\in H^2(\Omega)\cap V_\sigma$ and $p\in V$ such that $-\Delta u + \nabla p = f$ almost everywhere in $\Omega$.  Additionally,
\begin{equation}
    ||u||_{H^2(\Omega)} + ||p||_V \leq C ||f||
\end{equation}
for some $C>0$.

\subsection*{D.2 Proof}
\begin{proof}
First, we define the residual functions
\begin{equation}
    f_1 := \partial_t\tu + (\tu\cdot\nabla)\tu - \di(\nu(\tphi)D\tu) + \nabla P_\theta - \tmu\nabla\tphi
\end{equation}
\begin{equation} 
    f_2 := \partial_t\tphi + \tu\nabla\tphi - \Delta\tmu .
\end{equation}
Then we define the error or difference functions
\begin{equation}
    \bu := u - \tu \text{ , } \bphi = \phi - \tphi \text{ , } \bar{P} = P - P_\theta.
\end{equation}
The following equations are then satisfied
\begin{subequations} \label{NSCH:Difference:MainEq}
    \begin{numcases}{}
        \begin{aligned}
            \bu_t + (u\cdot\nabla)\bu &+ (u\cdot\nabla)\tu  - \text{div}((\nu(\phi)- \nu(\tphi))D\tu) \\ &- \di(\nu(\phi)D\bu)
            + \nabla \bar{P} = \mu\nabla\phi - \tmu\nabla\tphi - f_1
        \end{aligned} \label{NSCH:diff:ueq}\\
        \nabla\cdot \bu = 0 \\
        \bphi_t + (u\cdot\nabla)\bphi + \bu\cdot\nabla\tphi - \Delta\bmu = f_2  \label{NSCH:diff:phieq}\\
        \mu = -\Delta\phi + \Psi'(\phi)
    \end{numcases}
\end{subequations}
with Initial and Boundary Conditions

\begin{subequations} 
    \begin{numcases}{}
        \bu(\cdot,0) = u_0 - \tu(\cdot,0)  \text{ , } \bphi(\cdot, 0) = \phi_0 - \tphi(\cdot,0)  \text{ in } \Omega  \\
        \bu=0 \text{ , } \partial_n\bphi = 0 \text{ , }  \partial_n\bmu = 0 \text{ on } \partial\Omega\times (0,T) \label{NSCH:bmu:boundary}.
    \end{numcases}
\end{subequations}

Notice that although we express the above in the notation of strong solutions, we immediately pass to a weak formulation for our analysis, and in fact only a weak solution as defined in the preliminaries is necessary.  We obtain a priori estimates on these equations.  First, multiply equation (\ref{NSCH:diff:phieq}) by $A_0^{-1}(\bphi - \lbphi)$ where $\lbphi = |\Omega|^{-1}(\bphi, 1)$ (i.e. the average value of $\bphi$) and $A_0^{-1}$ is the inverse of the Neumann Operator described in the preliminaries.  From this we see
\begin{align} \begin{split}
    (\partial_t\bphi, \AZ) + (u\cdot\bphi,\AZ) + (\bu\cdot\nabla\tphi,\AZ) \\+ (\nabla\bmu,\nabla\AZ = (f_2,\AZ).
\end{split} \end{align}
Then rewrite to
\begin{equation} \label{NSCH:pre_sum_1} \begin{split}
    \langle\partial_t\bphi , \AZ\rangle - (\bphi u , \nabla\AZ) - (\tphi \bu , \nabla \AZ)\\ + (\nabla\bmu,\nabla\AZ) = (f_2, \AZ).
\end{split} \end{equation}
Now we analyze term by term.  First,
\begin{align} \begin{split}
    \langle\partial_t\bphi , \AZ\rangle &= \langle\partial_t(\bphi - \lbphi),\AZ\rangle + \langle \partial_t\bphi,\AZ\rangle \\
    &= \frac{1}{2}\frac{d}{dt}||\bphi-\lbphi||_*^2 + \langle\partial_t\bphi_{avg},\AZ\rangle \\
\end{split} \end{align}
where we recall that $||\bphi - \lbphi||_* = ||\nabla \AZ||_{L^2(\Omega)}$. We handle the second term later.  We label the next terms
\begin{equation}
    I_1 := (\bphi u , \nabla\AZ)
\end{equation}
\begin{equation}
    I_2 := (\tphi \bu , \nabla \AZ).
\end{equation}
Then treat the fourth term,
\begin{equation}
    (\nabla\bmu, \nabla\AZ) = -\int_\Omega \bmu \Delta \AZ dx + \int_{\partial\Omega}\nabla\AZ \bmu\cdot n d\Gamma
\end{equation}
by properties of $A_0^{-1}$ described in the preliminaries as well as the boundary conditions of $\bmu$, (\ref{NSCH:bmu:boundary}).  Then, 
\begin{equation} \label{NSCH:IM1}
    (\bmu , \bphi - \lbphi) = (-\Delta\bphi, \bphi - \lbphi) + (\Psi'(\phi) - \Psi'(\tphi),\bphi - \lbphi)
\end{equation}
The first term on the right side of (\ref{NSCH:IM1}) becomes
\begin{align} \begin{split}
    (-\Delta\bphi,\bphi - \lbphi) &= \int_\Omega (\nabla(\bphi-\lbphi))^2 dx - \int_{\partial\Omega} \frac{\partial (\bphi - \lbphi)}{\partial n} (\bphi - \lbphi) \\
    &= ||\nabla\bphi||_{L^2(\Omega)}^2
\end{split} \end{align}
by (\ref{NSCH:bmu:boundary}) and the fact that $\lbphi$ is a constant in x.  The second term of on the right side of (\ref{NSCH:IM1}) we bound using (\ref{PF:alpha0}-\ref{PF:alpha}) and the Fundamental Theorem of Calculus:
\begin{align} \begin{split}
    (\Psi'(\phi) - \Psi'(\tphi),\bphi - \lbphi) &\geq -\alpha (\bphi, \bphi - \lbphi) \\
    &=  - \alpha ||\bphi - \lbphi|| + \alpha (\bphi_{avg}, \bphi - \lbphi),
\end{split} \end{align}
but $(\bphi_{avg}, \bphi - \lbphi) = \int_\Omega \bphi\lbphi - |\lbphi|^2 dx = |\lbphi|^2(|\Omega| - |\Omega|) = 0$.

From (\ref{NSCH:pre_sum_1}), all we are left to handle is $( f_2,\AZ) - \langle \partial_t \lbphi , \AZ \rangle $.  By Integrating (\ref{NSCH:diff:phieq}), we see
\begin{equation}
    \int_\Omega \partial_t\bphi dx+ \int_\Omega u\cdot\nabla\bphi dx+ \int_\Omega \bu\cdot\nabla\tphi dx= \int_\Omega f_2 dx
\end{equation}
by noting that $\Delta\mu$ integrates to 0 by integration by parts with the boundary conditions.  However,
\begin{equation}
    \int_\Omega u\cdot\nabla\bphi = \int_{\partial\Omega} u \bphi\cdot n - \int_\Omega (\nabla\cdot u) = 0
\end{equation}
by the Divergence-free condition and the boundary condition.  Similarly,
\begin{equation}
    \int_\Omega \bphi\cdot\nabla\tphi = 0.
\end{equation}
So see
\begin{equation}
    \partial_t\bphi_{avg} = f_{2,avg}.
\end{equation}
Then if we label
\begin{equation}
    I_3 = (f_2 - f_{2,avg} , \AZ ),
\end{equation}
we get the following intermediary result
\begin{equation} \label{NSCH:sum_1}
    \frac{1}{2}\frac{d}{dt}||\bphi-\lbphi||_* + \frac{1}{2}||\nabla\bphi||^2 \leq \frac{\alpha^2}{2}||\bphi - \lbphi||_* + I_1 + I_2 + I_3.
\end{equation}
Now we aim towards a second intermediary result.  We take the inner product in $L^2$ of (\ref{NSCH:diff:ueq}) with $A^{-1}u$ as defined in the preliminaries.  Using the boundary conditions we get that
\begin{equation} \label{NSCH:sum_2}
    \frac{1}{2}\frac{d}{dt} ||u||_\#^2 + (\nu(\phi)\D{\bu}, \nabla A^{-1}u) = I_4 + I_5 + I_6 + I_7 + I_8
\end{equation}
where
\begin{align}  \begin{split}
    I_4 &= -((\nu(\phi) - \nu(\tphi))\D\tu, \nabla A^{-1}\bu), \\
    I_5 &= (u\otimes \bu , \nabla A^{-1}\bu) + (u\otimes\tu, \nabla A^{-1}\bu), \\
    I_6 &= (\nabla\phi\otimes\nabla\bphi , \nabla A^{-1}\bu) + (\nabla\bphi\otimes\nabla\tphi,\nabla A^{-1}\bu), \\
    I_7 &= (f_1. A^{-1}\bu), \\
    I_8 &= -(\nabla \bar{P}, A^{-1}\bu).
\end{split}  \end{align}

Then by observing
\begin{align} \label{NSCH:v_phi_eq_1} \begin{split}
    (\nu(\phi)\D \bu, \nabla A^{-1}\bu) &= (\nabla u, \nu(\phi)\D A^{-1}\bu) \\
    &= -(\bu, \di (\nu(\phi)\D A^{-1}\bu)) \\
    &= - (\bu, \nu'(\phi)\D A^{-1}\bu\nabla\phi) - \frac{1}{2}(\bu, \nu(\phi)\Delta A^{-1}u).
\end{split} \end{align}

We know from Stokes Operator Theory (see D.1.4) that $\exists p\in L^2(0,T;V)$ such that $-\Delta A^{-1}\bu + \nabla p = \bu$ where

\begin{equation}
    ||p||\leq c||\nabla A^{-1}\bu||^\frac{1}{2}||p||_V\leq c||\bu||.
\end{equation}
Thus
\begin{align} \label{NSCH:v_phi_eq_2}
\begin{split}
    -\frac{1}{2} (\bu,\nu(\phi)\Delta A^{-1}\bu) &= \frac{1}{2}(\nu(\phi)\bu),\bu) - \frac{1}{2}(\nu(\phi)\bu,\nabla p) \\
    &\geq v_*||u||^2 + \frac{1}{2}(\nu'(\phi)\nabla \phi\cdot\bu , p).
\end{split}
\end{align}

Combining (\ref{NSCH:v_phi_eq_1}) and (\ref{NSCH:v_phi_eq_2}),
\begin{equation} \label{NSCH:sum_3}
    (\nu(\phi)\D\bu, \nabla A^{-1}\bu) \geq v_* ||\bu||^2 - I_9 - I_{10}
\end{equation}
where
\begin{equation}
\begin{split}
    I_9 = (\bu , \nu'(\phi)\D A^{-1}\bu\nabla\phi) \\
    I_{10} = -\frac{1}{2}(\nu'(\phi)\nabla\phi\cdot\bu, p).
    \end{split}
\end{equation}
Then combing (\ref{NSCH:sum_1}), (\ref{NSCH:sum_2}), and (\ref{NSCH:sum_3}) we get the primary differential inequality
\begin{equation}
    \frac{d}{dt}H + \nu_*||\bu||^2 + \frac{1}{2}||\nabla\bphi||^2 \leq \alpha^2 H + \sum_{k=1}^{10} I_k
\end{equation}
where
\begin{equation}
    H(t) = \frac{1}{2}||\bu(t)||^2_\# + ||\bphi - \lbphi||^2_*.
\end{equation}
Now we attempt to control the $I_k$ terms.  By Holder's inequality, Sobolev Embedding, Poincaré's Inequality, and Cauchy's Inequality, we can bound $I_1$
\begin{align}
\begin{split}
        I_1 &= (\bphi u, \nabla \AZ) \\
        &\leq ||\bphi||_{L^6} ||u||_{L^3} ||\bphi-\lbphi||_* \\
        &\leq C_1 ||\nabla\phi ||_{L^2} ||u||_{L^3} ||\bphi-\lbphi||_* \\
        &\leq \frac{1}{8}||\nabla\bphi||^2 + C_1 ||u||^2_{L^3}||\bphi-\lbphi||^2_*.
\end{split}
\end{align}
We use Holder and Cauchy's inequality to bound $I_2$
\begin{align}
\begin{split}
    I_2 &= (\tphi\bu, \AZ) \\
    &\leq ||\tphi||_{L^\infty}||\bu||||\bphi-\lbphi||_* \\
    &\leq \frac{v_*}{8}||\bu||^2 + C_2 ||\bphi - \lbphi||^2_*.
\end{split}
\end{align}
Now for $I_3$,
\begin{align}
    \begin{split}
        I_3 &= \langle f_2-f_{2,avg} , \AZ \rangle \\
        &\leq \int_\Omega \nabla A_0^{-1} (f_2 - f_{2,avg}) \nabla A_0^{-1}(\bphi - \lbphi) dx\\
        &\leq ||f_2 - f_{2,abg}||_*||\bphi - \lbphi||_* \\
        &\leq \frac{1}{2} ||f_2 - f_{2. avg}||^2 + C_3||\bphi - \lbphi||^2_*
    \end{split}
\end{align}
by the properties of $A_0^{-1}$ described in the preliminaries (D.1.4) and by the equivalence of $||\cdot||_*$ to $||\cdot||$.  We move next to $I_5$, leaving $I_4$ for last.
\begin{align}
    \begin{split}
        I_5 &= (u\otimes \bu, \nabla A^{-1}\bu ) + (\bu\otimes \tu, \nabla A^{-1}u) \\
        &\leq (||u||_{L^4} + ||\tu||_{L^4}) ||\bu|| ||\nabla A^{-1}\bu||_{L^4} \\
        &\leq c(||u||^{\frac{1}{2}}||u||_V^{\frac{1}{2}}+ ||\tu||^{\frac{1}{2}}||\tu||_V^\frac{1}{2})||\bu|| ||\nabla A^{-1}\bu||_{L^4} \\
        &= c(||u||^{\frac{1}{2}}||u||_{V_\sigma}^{\frac{1}{2}}+ ||\tu||^{\frac{1}{2}}||\tu||_{V_\sigma}^\frac{1}{2})||\bu|| ||\nabla A^{-1}\bu||_{L^4} \\
        &\leq \frac{v_*}{8}||u||^2 + C_4(||u||^2_{V_\sigma} + ||\tu||^2_{V_\sigma})||\bu||^2_\#.
    \end{split}
\end{align}
Where we use Holder and Gagliardo-Nirenberg (\ref{Ineq:GN}) Inequalities, followed by the fact that $u$ and $\tu$ are divergence free.  We treat $I_6$ by
\begin{align}
    \begin{split}
        I_6 &= (\nabla\phi\otimes\nabla\bphi, \nabla A^{-1}\bu) + (\nabla\bphi\otimes\tphi, \nabla A^{-1}\bu) \\
        &\leq (||\nabla\phi||_{L^\infty} + ||\nabla\tphi||_{L^\infty})||\bphi||||\nabla A^{-1}u|| \\
        &\leq \frac{1}{8} ||\nabla\bphi||^2 + C_5(||\nabla \phi||^2_{L^\infty} + ||\nabla \tphi||^2_{L^\infty})||u||_\#^2.
    \end{split}
\end{align}
Then by Poincaré and Cauchy Inequalities
\begin{align}
    \begin{split}
        I_7 &= (f_1, A^{-1}\bu) \\
        &\leq c||f_1||||\nabla A^{-1}\bu|| \\
        &\leq c ||f_1|| ||\bu||_\# \\
        &\leq \frac{||f_1||^2}{2} + C_5 ||u||^2_\#.
    \end{split}
\end{align}
Since the spaces are perpendicular, 
\begin{equation}
    I_8 = -(\nabla P, A^{-1}\bu ) = 0.
\end{equation}
Then,
\begin{align}
    \begin{split}
        I_9 &= -\frac{1}{2}(\nu'(\phi)\nabla\phi\cdot\bu, p) \\
        &\leq c||\bu|| ||\D A^{-1}\bu|| ||\nabla\phi||_{L^\infty} \\
        &\leq \frac{\nu_*}{8} ||\bu||^2 + C_6 ||\nabla \phi||^2_{L^\infty}||u||_\#^2.
    \end{split}
\end{align}

For $I_{10}$, we use Agmon's Inequality (\ref{Ineq:AG}), Gagliardo-Nirenberg Inequality (\ref{Ineq:GN}), the bound of $p$ (\ref{Ineq:pbound}), and Cauchy Inequality,
\begin{align}
    \begin{split}
        I_{10} &= -\frac{1}{2}(\nu'(\phi)\nabla\phi\cdot u, p) \\
        &\leq c ||\nabla\phi||_{L^4} ||\phi||^{\frac{1}{2}} ||\bu|| ||p||_{L^4} \\
        &\leq c ||\phi||^\frac{1}{2}_{L^\infty} ||\phi||^\frac{1}{2}_{H^2} ||\bu|| ||p||_{L^4} \\
        &\leq c||\phi||^\frac{1}{2}_{H^2}||\bu|| ||p||^\frac{1}{2} ||p||^\frac{1}{2}_V \\
        &\leq c||\phi||_{H^2}^\frac{1}{2}||\bu||||\nabla A^{-1}\bu||^\frac{1}{2} ||\bu||^\frac{1}{4} ||\bu||^\frac{1}{2} \\
        &= c ||\phi||_{H^2}^\frac{1}{2} ||\bu||^\frac{7}{4} ||\nabla A^{-1}\bu || \\
        &\leq \frac{\nu_*}{8} ||\bu||^2 + C_7 ||\phi||^4_{H^2}||u||^2_\#.
    \end{split}
\end{align}

Finally, we treat $I_4$.
\begin{align}
    \begin{split}
        I_4 &= ((\nu(\phi)-\nu(\tphi))\D \tu, \nabla A^{-1}\bu ) \\
        &= \left(\int_0^1 \nu'(s\phi + (1-s)\tphi) ds \bphi\D\tu, \nabla A^{-1}\bu\right) \\
        &\leq c ||\D \tu|| ||\bphi \nabla A^{-1}\bu|| \\
        &\leq c||\D\tu|| ||\nabla \bphi||(||\nabla A^{-1}\bu|| + ||\bphi||)\left[ \log\left(c \frac{||\bu|| + ||\bphi||_V}{||\nabla A^{-1}\bu|| + ||\bphi||}     \right)\right]^\frac{1}{2} \\
        &\leq \frac{1}{4}||\nabla\bphi||^2 + C_8 ||\D \tu||^2 H \log\left(\frac{C_9}{H}\right).
    \end{split}
\end{align}
We note $H$ is bounded by the regularity theorem (\ref{AppD:Reg}), so we choose $C_9$ so that $\log\left(\frac{C_9}{H}\right)\geq 1$.  Combining the estimates for the $I_k$, we arrive at
\begin{equation}
    \frac{d}{dt}H \leq H Y \log\left(\frac{C_9}{H}\right) + \frac{1}{2} ||f_2 - f_{2, avg}||^2 + ||f_1||^2,
\end{equation}
where
\begin{equation}
    Y(t) = C_9\left( 1 + ||u||^2_{L^3} + ||u||^2_{v_\sigma} + ||\nabla\phi||^2_{L^\infty} + ||\nabla\tphi||^2_{L^\infty} + ||\phi||^4_{H^2}\right).
\end{equation}

So by the boundedness of the Loss function
\begin{equation}
    H'(t) \leq H(t)Y(t)\log\left(\frac{C_9}{H}\right) + \frac{3}{2}\epsilon^2.
\end{equation}
We multiple both sides of this equation by $\exp(-\int_0^t Y(\tau)\log C_9/H(\tau)) d\tau)$.  Then,
\begin{align} \begin{split}
    H'(t)e^{-\int_0^t Y(\tau)\log C_9/H(\tau)) d\tau} - H(t)Y(t)\log(C_9/H(t))e^{-\int_0^t Y(\tau)\log(C_9/H(\tau))d\tau} \\ \leq \frac{3C_{10}}{2}\epsilon^2e^{-\int_0^t Y(\tau)\log C_9/H(\tau) d\tau},
\end{split} \end{align}
which is equivalent to 
\begin{equation}
    \left[ H(t)e^{-\int_0^t Y(\tau)\log C_9/H(\tau)) d\tau}\right]' \leq \frac{3}{2}\epsilon^2e^{-\int_0^t Y(\tau)\log C_9/H(\tau) d\tau}.
\end{equation}
We then infer
\begin{align}
\begin{split}
    H(t)e^{-\int_0^t Y(\tau)\log C_9/H(\tau) d\tau} &\leq H(0) + \frac{3}{2}\epsilon^2\int_0^te^{-\int_0^s Y(\tau)\log C_9/H(\tau) d\tau}ds \\
    &\leq \epsilon^2(1+\frac{3}{2}t).
\end{split}
\end{align}
Since $\log(C_9/H) > 0$ and $Y(t) > 0$,
\begin{equation}
    H(t)e^{\int_0^t Y(\tau)\log H(\tau) d\tau} \leq \epsilon^2(1+\frac{3}{2}t)e^{\log(C_9)\int_0^tY(\tau)d\tau}.
\end{equation}
Finally, if we let $\epsilon \to 0$, we see that $H(t) \to 0$ since the exponential is not 0.
\end{proof}

\pagebreak

%\addcontentsline{toc}{section}{References}
%\nocite{*}
%\bibliography{annot} 

\begin{thebibliography}{99}
\bibitem{heattransfer}S. Cai, Z., Wang, S. Wang, P. Perdikaris, and G. Karniadakis, {\it Physics-Informed Neural Networks for Heat Transfer Problems}, Journal Of Heat Transfer. \textbf{143}, 060801 (2021, 4).

\bibitem{highspeed}Z. Mao, A. Jagtap, and G. Karniadakis, {\it Physics-informed neural networks for high-speed flows}. Computer Methods In Applied Mechanics And Engineering \textbf{360} pp. 112789 (2020).

\bibitem{cardiac}F. Sahli Costabal, T. Yang, P. Perdikaris, D. Hurtado and E. Kuhl, {\it Physics-Informed Neural Networks for Cardiac Activation Mapping}, Frontiers In Physics \textbf{8} (2020).

\bibitem{radiative}S. Mishra and R. Molinaro, {\it Physics informed neural networks for simulating radiative transfer}, Journal Of Quantitative Spectroscopy And Radiative Transfer \textbf{270} (2021).

\bibitem{poroelasticity}E. Haghighat, D. Amini and R. Juanes, {\it Physics-informed neural network simulation of multiphase poroelasticity using stress-split sequential training}, Computer Methods In Applied Mechanics And Engineering \textbf{397} pp. 115141 (2022).

\bibitem{atomistic}G. Pun, R. Batra, R. Ramprasad and Y. Mishin, {\it Physically informed artificial neural networks for atomistic modeling of materials}, Nature Communications \textbf{10}, 2339 (2019,5).

\bibitem{fasttime}J. Stiasny, S. Chevalier and S. Chatzivasileiadis, {\it Learning without Data: Physics-Informed Neural Networks for Fast Time-Domain Simulation.} 2021 IEEE International Conference On Communications, Control, And Computing Technologies For Smart Grids (SmartGridComm) pp. 438-443 (2021)

\bibitem{crustal}T. Okazaki, T, Ito, K. Hirahara and N. Ueda, {\it Physics-informed deep learning approach for modeling crustal deformation}, Nature Communications \textbf{13}, 7092 (2022,11)

\bibitem{fluids8020046}P. Moser, W. Fenz, S. Thumfart, I. Ganitzer and M. Giretzlehner, {\it Modeling of 3D Blood Flows with Physics-Informed Neural Networks: Comparison of Network Architectures}, Fluids \textbf{8} (2023).

\bibitem{HAL}N. Doumèche, G. Biau and C. Boyer, {\it Convergence and error analysis of PINNs}, Hal-04085519v1 (2023).

\bibitem{residuals}Y. Shin, Z. Zhang and G. Karniadakis, {\it Error Estimates of Residual Minimization using NNs for Linear PDEs}, Journal Of Machine Learning For Modeling And Computing \textbf{4} (2023,1).

\bibitem{Convergence}J. Shin and G. Karniadakis, {\it On the Convergence of Physics Informed Neural Networks for Linear Second-Order Elliptic and Parabolic Type PDEs}, Communications In Computational Physics \textbf{28}, 2042-2074 (2020).

\bibitem{RAISSI2019686}M. Raissi, P. Perdikaris and G. Karniadakis, {\it Physics-informed neural networks: A deep learning framework for solving forward and inverse problems involving nonlinear partial differential equations},  Journal Of Computational Physics \textbf{378} pp. 686-707 (2019).

\bibitem{Liu1989}D. Liu and J. Nocedal, {\it On the limited memory BFGS method for large scale optimization}, Mathematical Programming \textbf{45}, 503-528 (1989,8).

\bibitem{adam}D. Kingma and J. Ba, {\it Adam: A Method for Stochastic Optimization}, International Conference On Learning Representations (2014,12).

\bibitem{autodiff}A. Baydin, B. Pearlmutter, A. Radul and J. Siskind, {\it Automatic differentiation in machine learning: a survey}, J. Mach. Learn. Res. \textbf{18}, 5595-5637 (2017,1).

\bibitem{krishnapriyan2021characterizing}A. Krishnapriyan, A. Gholami, S. Zhe, R. Kirby and M. Mahoney, {\it Characterizing possible failure modes in physics-informed neural networks} (2021)

\bibitem{Biswas2022}A. Biswas, J. Tian and S. Ulusoy, {\it Error estimates for deep learning methods in fluid dynamics}, Numerische Mathematik \textbf{151}, 753-777 (2022,7).

\bibitem{wight2020solving}C. Wight, and J. Zhao, {\it Solving Allen-Cahn and Cahn-Hilliard Equations using the Adaptive Physics Informed Neural Networks}, Communications in Computational Physics  (2020).

\bibitem{haitsiukevich2023improved}K. Haitsiukevich and A. Ilin, {\it Improved Training of Physics-Informed Neural Networks with Model Ensembles}, International Joint Conference on Neural Networks (2023).

\bibitem{MATTEY2022114474}R. Mattey and S. Ghosh, {\it A novel sequential method to train physics informed neural networks for Allen Cahn and Cahn Hilliard equations}, Computer Methods In Applied Mechanics And Engineering, \textbf{390} pp. 114474 (2022).

\bibitem{zhang2024priori}G. Zhang, J, Lin, Q. Zhai, H. Yang, X. Chen, X. Zheng and I. Leong, {\it A Priori Error Estimation of Physics-Informed Neural Networks Solving Allen–Cahn and Cahn–Hilliard Equations} ArXiv (2024).

\bibitem{Sun2019SurrogateMF}L. Sun, H. Gao, S. Pan and J. Wang, {\it Surrogate modeling for fluid flows based on physics-constrained deep learning without simulation data}, Computer Methods In Applied Mechanics And Engineering \textbf{361} (2019).

\bibitem{chudomelka2020deep}B. Chudomelka, Y. Hong, H. Kim and J. Park, {\it Deep neural network for solving differential equations motivated by Legendre-Galerkin approximation}, International Journal of Numerical Analysis and Modeling \textbf{21} (2023).

\bibitem{ko2022convergence}S. Ko, S. Yun and Y. Hong, {\it Convergence analysis of unsupervised Legendre-Galerkin neural networks for linear second-order elliptic PDEs} CoRR (2022).

\bibitem{RichterPowell2022NeuralCL}J. Richter-Powell, Y. Lipman and R. Chen, {\it Neural Conservation Laws: A Divergence-Free Perspective}, ArXiv (2022).

\bibitem{deryck2023error}T. Ryck, A. Jagtap and S. Mishra, {\it Error estimates for physics informed neural networks approximating the Navier-Stokes equations}, IMA Journal of Numerical Analysis \textbf{44} (2023).

\bibitem{Jin_2021}X. Jin, S. Cai, H. Li and G. Karniadakis, {\it NSFnets (Navier-Stokes flow nets): Physics-informed neural networks for the incompressible Navier-Stokes equations}, Journal Of Computational Physics, \textbf{426} pp. 109951 (2021,2)

\bibitem{Cuomo2022}S. Cuomo, V. Di Cola, F. Giampaolo, G. Rozza, M. Raissi and F. Piccialli, {\it Scientific Machine Learning Through Physics–Informed Neural Networks: Where we are and What's Next}, Journal Of Scientific Computing \textbf{92}, 88 (2022,7).

\bibitem{doi:10.1137/18M1223459}A. Giorgini, A. Miranville and R. Temam, {\it Uniqueness and Regularity for the Navier–Stokes–Cahn–Hilliard System}, SIAM Journal On Mathematical Analysis \textbf{51}, 2535-2574 (2019).

\bibitem{GIORGINI2020194}A. Giorgini, and R. Temam, {\it Weak and Strong Solutions to the Nonhomogeneous Incompressible Navier-Stokes-Cahn-Hilliard system}, Journal De Mathématiques Pures Et Appliquées \textbf{144} pp. 194-249 (2020).

\bibitem{NSCHattractors}A. Giorgini, and R. Temam, {\it Attractors for the Navier-Stokes-Cahn-Hilliard system}, Discrete And Continuous Dynamical Systems - S \textbf{15}, 2249-2274 (2022).

\bibitem{Chen2022}W. Chen, J. Jing, C. Wang, and X. Wang, {\it A Positivity Preserving, Energy Stable Finite Difference Scheme for the Flory-Huggins-Cahn-Hilliard-Navier-Stokes System}, Journal Of Scientific Computing, \textbf{92}, 31 (2022,6).


\end{thebibliography}
%\bibliographystyle{plain}

\end{document}